**Werner DePauli-Schimanovich[1]**
**Paradigms-Shift in Set Theory, Ch64-SE9**

(0) Abstract
In this article the author claims that there is a paradigm shift from ZFC to NFUM and further to NACT - due to philosophical reasons, not mathematical ones. The goal is to construct systems where every "Not-Properclass" is a set! With help of Non-Monotonic Logic, the consistent systems NACT-MoonW, NACT*W, and NACT-SunW are producing "largest possible universes" of sets. Using self-evident philosophical principles, three approximations are suggested: NACT+NFUM-closed (to NACT*), NACT&ZFC4+(GCH) (also to NACT*) and NACT-NFUM (to NACT-Sun). Also the system NACT[&ZFC4-closed]+(FCA) is considered.

NFUM = NFU with (AC) and measurable properclass Ord. NFUM-closed is NFUM where the set-constituting formulas $A_i$ in the set operator need not only to be stratified but also to be made parameter-free (i.e., have only x as free variable over which the set is comprehended). In other formulas, free variables are allowed.

Keywords: Set Theory, ZFC, Naïve Set Theory, Predicate-Extension, Church Schema of Comprehension, NF, NFUM, Stratification, Universal Sets, Eliminative Class Theory, NBG, NACT (= Naïve Axiomatic Class Theory), NACT-Sun, NACT-Moon, NACT-Star, NACT+NFUM-closed, ZFCK, NACT+NFUM, (GCH), (FCA) [= Finite or Countable Anti-Thesis], NACT+ZFC4+(GCH), NACT[+ZFC4-[closed]+(FCA).

(1) Introduction
When mathematicians hear the word "Set Theory", they think first of the axiom system of Zermelo-Fraenkel ZFC [i.e. ZF with the axiom of choice (AC)], despite the fact that there probably exist 100 similar systems.

After the contradictions of Naïve Set Theory, only the most important principles of set-construction necessary for the construction of mathematics have been accepted. These are: the axiom of extensionality, the existence of the empty set, of pairs and unions, of an infinite set, of parts of sets, of the power set, the large union ( = sum) of a set, and of images of a set.[2] This fixing of the axioms of ZFC is of course somewhat arbitrary. But David Hilbert[3] wanted to show (with his program of Finitism) first the consistency of Peano-Arithmetics and then of ZFC. Since Gödel[4] we know that this program has been toppled, and an arbitrary limitation of the axioms is no longer justified.

---

[1] Department for Data Bases & A.I., Institute for Information Systems, Technical University Vienna, Austria, Werner.DePauli@gmail.com

[2] Some of these axioms are redundant, because it can be demonstrated that they follow logically from the other ones.
[3] See Richard Zach in the literature.
[4] See my books in the literature, Smullyan [1992], Smorynski [2002] and Goldstern [1998].

When philosophers hear the word "Set Theory", they think first of a general principle of set-construction, since they want to construct the world "from a logical point of view".[5] The Naïve Set Theory was just such a general system of set-construction. But Bertrand Russell's Theory of Types was also such a system, from which Willard van Orman Quine has developed his system "New Foundations" NF. From this it can be seen that the goal of Logicism is to constitute a set theory with nearly "unrestricted" Church Schema of set-comprehension.

The all-decisive keyword for a philosophically-based set theory is the "Predicate-Extension": $y \in \{x: A(x)\} \iff A(y)$. As we know, this principle was originally considered unlimitedly valid for sets too, and made it possible for Russell to construct his famous paradox with the Russell-set $ru == \{x: x \notin x\}$.[6]

After that, the unrestricted predicate extensions were renamed and called "Classes" from then on. The elements of the classes can only be some special classes now called "Sets" (and which are of course also predicate extensions). This was the beginning of Neumann-Bernays-Gödel Class Theory NBG. For the fixation of sets, the arbitrarily selected axioms of ZFC are still used.

But classes are (in most cases) only virtual entities, which can even be eliminated with the help of the Church schema (in some class systems). Their meaning is nothing other than the predicate extension. The objects that truly exist are the sets! To distinguish the class operator from the set operator, we put a vertical bar after the opening curly set bracket (i.e. brace) and one before the closing one. So $\{|x: A(x)|\}$ produces always a class, but $\{x: A(x)\}$ produces a set only for appropriate A's or the empty set.[7]

(2) NF and NFUM

Therefore, if we are also looking for new sets (not included in ZFC) we require set theories with an only slightly restricted Church schema. The first one confronting us with that is New Foundation. [Unfortunately NF is very complicated and therefore only used by a few mathematical logicians. Furthermore, it is incompatible with the (AC).] Here the Church schema for sets is valid only for "Stratified" predicates $A_i$. Such an only slightly restricted comprehension schema entails of course the existence of very many sets (similar to Naïve Set Theory). Most of them are not required for mathematics; but they are very important for philosophy. For example, NF has a universal set $v == \{x: x = x\}$ and many sets equipollent (i.e. circa with the same size) to it.

---

[5] See Quine's book of the same name in the literature.

[6] There also exist methods with special restricted Church schemata to convert these properclasses into sets too. But these methods require an unjustified technical effort, which is the reason why we do not want to consider them here. Still, see also: "On Frege's True Way Out" in DePauli [1990] and "Naïve Axiomatic Mengenlehre for Experiments" in DePauli [2006e].

[7] This is a new convention suggested by the author.

Preferable to NF is the system NFUM of Randall Holmes[8], which is NF enlarged for "Ur-Elements" (or "Atoms"). These entities were first introduced to NF in the famous article "On the consistency of a slight (?) modification of Quine's New Foundation" by Ronald Björn Jensen.[9] As opposed to NF, NFU is compatible with (AC). NFUM contains (AC) and has also some further axioms and principles that make it easy to develop mathematics similar to ZFC. For large collections (i.e. classes) of sets, Holmes shows step by step that the stratification condition can be dropped and set-comprehension can therefore be executed (for these classes) just as in Naïve Set Theory.

In collaboration with Robert Solovay, Randall Holmes constructs an inner model of NFUM in which ZFC can also be implemented.[10] The strength of this system is comparable to ZFC along with very strong cardinal axioms (e.g. that for every n a n-Mahlo cardinal exists). The M on the end of NFUM should remind us of the fact that the universe of all ordinals is "Measurable" here.

With NFUM, large areas of mathematics can be shown in the same way as with ZFC. But the axioms and concepts are defined in an adequate philosophical way. E.g., the number "1" is the set of all sets with exactly one element, and so on (as in the Theory of Types.)

(3) Comparison of ZFC and NFUM

In ZFC, the Church Schema (ChS) is generally valid ONLY for classes:

$y \in \{|x: A(x)|\} \Longleftrightarrow A(y)$ [for arbitrary A].

For a given $A_i$ it is only certain that its predicate extension is a class. [We wish to indicate this with the class operator $\{|x: \ldots|\}$.] Some of these classes that should also be a set are listed in a separate list of axioms. [We can also define the sets as the elements of classes. But from this definition, the existence of any set cannot be derived, not even of the empty set.]

ZFC generates (by its axioms) only the smallest common denominator of sets, respectively those collections which can be considered consistently as sets, e.g. sets that are generated by one of the approximately 100 systems of set theory [such as Ackermann's system, for example].[11] Therefore, I suggest calling such set theories in the future "minimal set theories".

The (ChS) of ZFC is also generally valid for the classes of NFUM (e.g. the Russell properclass). But in NF and NFUM it is also valid for sets if we restrict it to stratified predicates [i.e. such predicates where the variables can be indexed with natural numbers such that for every pair of variables left and right of an element-sign "$\in$" the index on the right side is 1 higher than the one on the left side. E.g., "$x_7 \in y_8$ & . . . & $z_6 \in x_7$" is correct, but "$x \in$

---

[8] See the book of Randall Holmes (for NFUM) in the literature.

[9] See the article of Jensen in the literature.

[10] With regard to this, Holmes writes on page 163 of his book: "A suitable redefinition . . . . will actually cause the membership relation E of our model to coincide with the 'real' membership epsilon in a corner of our universe which can then be taken to be the universe of ZFC."

[11] See Felgner [2002].

x" is not. It is evident that every variable-letter has to be assigned the same index-number.] Since this is a decidable property, we can use a set operator instead of the class operator, which always generates a set if the predicate is stratified [and 0 if not]. Let us call it {x: … } [without vertical bars] which yields (ChS-strat) for sets:

y ∈ {x: A(x)} <==> A(y)   [for stratified A].

   As opposed to ZFC, NFUM has with (ChS-strat) a unique set-comprehending system. Unfortunately it is a bit unnatural for mathematicians to use, and ten times as difficult as ZFC. Therefore, it has not yet achieved a breakthrough in practical use. But for philosophers, NFUM is more natural than ZFC, because the concepts of mathematics (sets, numbers, functions, etc) can be defined within NFUM in an adequate way. Furthermore, NFUM is a "maximal set theory" which produces "the greatest common multiple of sets" (so far known).

   (4) "Greatest Possible Maximal Set Theories"
The success of these (and other) investigations has encouraged the philosophical interest in "Set Theories with a Universal Set"[12] in recent years. We can pose the following provocative question: Does there exist a set theoretical system that produces the "most" sets? Our motto should be: "Every Not-Properclass shall be a set!" [This is of course only a heuristically formulated goal that will be subsequently developed. Not-Properclasses should be all classes for which it not necessarily (but only logically) entails that they are a properclass. E.g. the universal class is a properclass in ZFC but not in NF and many other systems.]

   Because of these reasons I have made the effort to construct set theoretical systems that also produce greatest common multiples of sets, but are nevertheless not as difficult as NFUM. These systems should be based on principles similar to Naïve Set Theory and should correspond to the mental procedure of mathematics. Furthermore, they should be useable in computational work. I call this group of systems NACT (= Naïve Axiomatic Class Theory). Of course it is not necessarily required here for these systems to be Formal Systems[13] in the sense of Hilbert-Gödel. In particular, they need not have a recursively enumerable set of axioms, i.e. they need not to be axiomatizable. Such systems exist e.g. in multivalued logic, logic with infinitely long expressions, etc. We have to define criteria for which formulas we are to consider as axioms and which not.

In the following, we wish to describe some systems of NACT. After explaining [the parameter-free systems] NACT ☻ W (= the naked moon Wittgenstein) and NACT ☻ (= NACT-Moon) we proceed to [the also parameter-free systems] NACT* (= NACT-Star) and NACT*W. After investigating the formal system NACT+NFUM-closed [which is an approximation of NACT*] we study NACT☼ (= NACT-Sun) [which allows also parameters] and its approximation NACT+NFUM. NACT+NFUM is the most important of these systems, because it is the greatest possible set theory yet known as a formal system. But we have to first explain the other systems, because it is easier for the reader to grasp the idea (or ideology) beyond these mechanisms.

---

[12] See the book of Thomas Forster (for Universal Sets and NF).
[13] See Smullyan [1961].

Before discussing NACT+NFUM, we shall also glance at the system NACT&ZFC4+(GCH), which is smaller than NACT+NFUM, but works with ZFC and without NFUM. After NACT+NFUM we shall also present the systems NACT [&ZFC4-closed] +(FCA), which is an alternative to the last mentioned systems (especially for physicists and other scientists) with only finite or countable sets.

(5) Naïve Axiomatic Class Theory

We consider the first system, NACT ☻ W, mainly for didactical reasons. It is not a Formal System because it has "undecided" axioms, i.e. formulas about which we do not know a priori if they are axioms or not and for which we do not (yet) have any decision procedure.[14] In practice, this turns out to function more or less as in Non-Monotonic Logic. [Therefore Gödel's Theorem cannot be applied here.] But this handicap (in my eyes) can be neglected, because NACT ☻ W is consistent.

The simplest consistent system is the "Wittgenstein-System": every time a contradiction appears, stop the derivation and cross out the contradiction![15] Therefore this system is per definition consistent (or better: inconsistency-free). NACT ☻ W works similar, but is not as simple. (The W shall remind us of the Wittgenstein-system.)

The main task of our NACT-systems is to develop for every system a kind of algorithm to discover its axioms. For all systems we have first to define a formula generator that produces the parameter-free wffs $A_i$ of set theory with the only free variable "x" (one wff after the other) corresponding to its length and depending on the additional number of bound variables. We call this algorithm GL1 (= Generation of List1):

Verum, Falsum,

$x \in x$, $x \notin x$,

$\exists x_1: x_1 \in x$, $\exists x_1: x_1 \notin x$, $\exists x_1: x \in x_1$, $\exists x_1: x \notin x_1$,

negation of the preceding formulas,

conjunction of 2 of them,

simplification and crossing out already produced formulas,

renaming of variables,

introduction of new quantified variables,

etc.

[The reader can see that this procedure is ideal for computational set theory.]

After some time the following formula will be produced:

---

[14] It is not impossible that such a decision procedure exists. But until now we have no hints for how to find it.

[15] Personal communication of Georg Kreisel 1983 to the author.

$\forall$ x₁: 0 $\in$ x₁ & $\forall$ x₂: x₂ $\in$ x₁ => x₂ $\cup$ {x₂} $\in$ x₁ ==> x $\in$ x₁. (This is the formula which generates omega.)
Let us call this sequence of wffs the original "list 1".

### (6) The System NACT ☻ W

In the following we substitute the formulas $A_i$ (one after the other) into the class operator {|x: $A_i$(x)|}, starting thus to form the original axiom-list1. After this, we assume that this class is a set, and try to derive a contradiction from every such assumption. If we succeed, we cross the assumed axiom out of the axiom-list1. Furthermore, we also cross out the assumed axiom with the negated formula ￢$A_i$ from axiom-list1, i.e. cross out Set({|x: ￢$A_i$(x)|}). [For reasons of computational efficiency, it would also be good to cancel $A_i$(x) and ￢$A_i$(x) from original list1.] We call now the lists checked lists until $A_i$.

This derivation of the contradiction is only allowed to be carried out within pure predicate logic FOL with identity and function-variables, in addition to an elementary set theory [but without any ZF-axioms except extensionality] that allows us to incorporate the already generated sets {x: $A_{i-j}$(x)} into the derivation [i.e. those classes {|x: $A_{i-j}$(x)|} for which the assumption that they are sets did not produce a contradiction].

Let us call this logic with class operator, Church Schema and axiom of extensionality the Class Theoretical Frame CT. If we include the (already mentioned) elementary set theory and the already derived sets into the proof of the "sethood" of a class we call this logic now CT&ES (= CT and Existing Sets). Our main task for every class is to derive falsum from the assumption that the class is a set:

Set({|x: $A_i$(x)|}) $\vdash_{\overline{CT\& ES}}$ Falsum.

If we succeed, $A_i$ and ￢$A_i$ will be cancelled from list1 and Set({|x: $A_i$(x)|}) and Set({|x: ￢$A_i$(x)|}) from axiom-list1. We can also install some decision criteria that a class is a set [e.g. for $A_i$(x) == x=x or for every complement of a set] and try to formalize some relative consistency proofs. Since we do not have unlimited time, the computation has to be done in time slices. But at the beginning, let us agree to consider a class in the axiom-list1 as a set as long as the inconsistency of this assumption was not shown (Non-Monotonic Reasoning).

Since all contradictory assumptions are removed immediately from the axiom-list1, no contradiction can be derived in NACT ☻ W and the system is therefore consistent. But we explained this system only for didactical reasons and do not want to use it in practice. The following simplifications will make the other systems without postfix W easier to work with.

Furthermore: I conjecture that in CT [e.g. NBG without ZFC] it is easy to show that a class is a properclass (like the Russell-class). In the article "The Notion 'Pathology' in Set Theory" (Chapter 10 and 11), the reader can find a more detailed investigation of these properclasses, called "pathologies" there.[16] If my conjecturing is true, it is possible to cross out all pathologies and their negations from the axiom-list1 and consider the rest of the axioms as valid. [This will be feasible especially in the next system NACT ☻ now to be investigated.]

---

[16] Definition of Hereditary-non-Patho :

(7) The System NACT ☻

We can simplify NACT ☻ W when we ignore the already proven sets. That yields NACT ☻ . Hereby we assume implicitly the philosophical thesis (J1) that establishing new sets is independent of already existing sets. [(J1) is the 3rd Patho-Thesis (3PT) of Chapter 9 of "The Notion 'Pathology' in Set Theory".] In NACT ☻ the main task is therefore simpler to prove, i.e. the derivation of falsum needs only be done within pure CT-Frame without ES (= considering the already proven sets). But therefore NACT ☻ is also stronger than NACT ☻ W and we do not know (any longer) if NACT ☻ is consistent.

NACT ☻ has also "undecided" axioms (as does NACT ☻ W), e.g. Set({|x: AC(x)|}), i.e. the class of sets which fulfil the axiom of choice. For all "pathological" classes it is of course easy to show the main task, e.g. Set({|x: x $\notin$ x|}) $\vdash_{CT}$ falsum. [To postulate axioms from which it is impossible to know if they are compatible with the other already "proven" axioms is not such a bad thing, and also occurs in ZFC: nearly all axioms concerning large cardinals are not known to be consistent. Here also the "wrongness" of an axiom can be shown more easily than its "validity" (for which we need a relative consistency proof).]

Instead of proving the main task for the whole list1, we can choose an arbitrary formula $A_i$ (for which we want to prove the sethood of the generated class) which is "consistent"[17], i.e. $A_i$ (and $\neg A_i$) are "hereditary not pathological". [I.e. if a wff $B_i$ or $\neg B_i$ produces a properclass, neither of them appear in $A_i$.] The collection of these $A_i$ is probably not recursively enumerable and therefore undecidable. [But it is not impossible that it is r.e.!]

NACT ☻ has no parameters in the set-constituting classes. [But in other formulas their use is allowed!] This is the reason why we cannot prove the theorem of Cantor. [I.e., that the power-set is higher in might than the original set.] Therefore we need the Generalized

---

(a) A wff is prim-patho (= primitive pathological) if it is a circle not-($x \in^n x$) or a Mirimanoff-formula M(x) [which says that x has an infinitely descending element-sequence]. My conjecture: these 2 types are the only primitive pathological wffs. And the set of prim-patho wffs is therefore decidable.
(b) A wff is patho(logical) if a prim-patho wff can be derived from it. [Research project: is this collection r.e.?]
(c) A wff is anti-patho if the negation of a prim-patho wff can be derived.
(d) A wff is hereditary patho if it contains a patho subformula.
(e) A wff is HnP (= hereditary-non-patho) if it contains no patho subformula. [If (b) is decidable, then also HnP.]
[17] Definition of "consistent":
(f) A wff is "consistent" if it contains neither a patho subformula nor an anti-patho subformula. ["Consistent" means here: the assumption that this wff (or its negation) generates a set does not lead to a contradiction.]

Continuum Hypothesis (GCH) to solve this situation. But another solution would be the system NACT ☺ (= NACT-Fullmoon), where no powers higher than countability exist. To establish this, we shall use a General Continuum Anti-Thesis (GCAnti-T1) [e.g. (FCA), the "finite or countable anti-thesis"]. This way [using (Anti-T2) or (Anti-T3), etc.] we can also define set theories with only finite sets or with finite, countable or continuous sets.[18] This would therefore yield a practical system for physicists or cosmologists who deny higher cardinalities.

(8)    The System NACT*

Because of the construction of list1, algorithm GL1 and axiom-list1 (where we take out primitive pathological $A_i$ and their negations), in NACT ☻ (and also in NACT ☻ W and NACT ☺ ) the complement (to every set) exists. If we now modify list1 and throw out only the "prim-patho" predicates (but not their negations), this results in new lists (to be called list2 and axiom-list2) and the new algorithm GL2. This is the system NACT*, which is a very important one.[19]

But this does not make the new solution of the main task "Set({|x: $A_i(x)$|}) $\vdash$ falsum" much easier (than before) because we have now more wffs to check. But this change is necessary. The reason for it is our doctrine "All Not-Properclasses should be sets!" And the complements of the properclasses can be considered as sets without any contradiction! In this case, of course, no general complements are allowed. Only small and large sets have general complements. But in any case, we cannot generate a general complement because no parameters are allowed in the class operators (which should be proven to generate a set).

As in NACT ☻ , we cannot prove Cantor's theorem in NACT* (because we have no free parameters in the class operator -- therefore, diagonalization can not be carried out and we have to fix this case again with (GCH)). With the help of an anti-thesis (GCAnti-H1) [like (FCA)] we can also define a variant of NACT* without higher cardinalities. Lets call it NACT ✴ (= NACT-Goldstar). Since in NACT* also the complements of the properclasses are sets, we use here implicitly the 1st Patho-Thesis (1PT) [or later (J3)]: For every pair of wffs A and ¬A, at least one of the two generates a set. This will help us to suggest an approximation to NACT* which is a Formal System. NACT* differs from NACT ☻ only in that we are now working with list2, GL2 and axiom-list2 (instead of 1).

(9)    NACT+NFUM() as Approximation to NACT*

Because of the undecidedness of some axioms (and the theoretical probability of undecidedability of hereditary-not-patho and "consistent" formulas) we can expect large computations in computational NACT*. To restrict these computations, we want to consider a Formal System with a r.e. axiom set which tries to save as many sets of NACT* as possible.

---

[18] See also the paper on "Pocket Set Theory" Holmes [2005].

[19] The reader can find a detailed discussion of NACT* in DePauli-Schimanovich [2006a]. NACT* (and later also NACT✫) are essentially Naïve Set Theory with Non-Monotonic Logic, where axioms are generated concerning its length.

We call it NACT+NFUM() and say "it converges" to NACT*. This new system is based on four philosophical principles, the truth of which I am conjecturing:

(J1)   The generation of new sets is independent of the already existing sets. [We used this principle by moving the W-systems to the one without W, i.e., to NACT ☻ and NACT*.]

(J2)   Small and large classes are sets. Small(X) ➔ Set(X) & Set($\neg$X). As small sets we use here the slim ones: slim(X) :== card(X) < card($\neg$X).

(J3)   Of 2 formulas $A_i$ and $\neg A_i$, at least one generates a set. I.e. Set({|x: $A_i$(x)|}) or Set({|x: $\neg A_i$(x)|}). [From this follows that the complements of the properclasses are sets.]

(J4)   The 2 principles (J2) and (J3) should be applied to NFUM(). [I.e., NFUM-closed[20] where the stratified formulas in the set operator of NFUM do not contain (free) parameters. In the other formulas of NFUM() parameters are of course allowed.] Stratified($A_i$()(x)) ➔ Set({|x: $A_i$()(x)|}). $A_i$ is a closed formula with the only free variable x. [The closure in the formula in the set operator (respectively class operator) is necessary to make it impossible to define the complement of a set!]

(J2) & (J3) & (J4) yields NACT+NFUM().

This recursively enumerable axiom-system produces universal sets and a lot of very large sets. The axiom of replacement (= image set axiom) and axiom of separation (= partial set axiom) are valid for small sets. [This follows directly from (J2).] Furthermore, the axiom of pairs is valid. The intersection and the union of small or large sets are again sets. About "medial" sets [with x ~ $\overline{x^1}$] not very much is known [except that the complements of the properclasses are medial]. Therefore the validity of laws for arbitrary intersections and unions is not clear. The situation for the axiom of sums is similar. [Medial sets are e.g. Slim := {|x: slim(x)|} or the classes of all singletons, pairs, triples, etc.]

In NACT+NFUM() we again cannot prove Cantor's theorem. So we have to assume (GCH) or some Anti-Thesis. We can also construct powersets only from the bottom as in the Cumulative Hierarchy because we do not allow parameters. Since we assumed in NFUM() implicitly the axiom of infinity, we can show the existence of $\omega_0$ and construct P($\omega_0$), etc. As opposed to NF, the (AC) is compatible with NFU and therefore also used as one of the axioms of NFUM. But it is not certain that (AC) is also compatible with NFUM(). Therefore it may be possible that only (AC-small) is valid [= (AC) only for small sets]. (The question of whether NF() is consistent is open to research. The same question arises now for NFUM(). NFU is relative consistent with ZF and therefore also NFU(). But this says nothing about the relative consistency of NFUM() with ZFC.)

As already mentioned in chapter 2, NFUM contains ZFC as an inner model. Therefore, all ZFC axioms are valid in NFUM, at least in a restricted form. But it is unknown if they are also valid in NFUM() and, for that matter, in NACT+NFUM(). It seems that they can be shown at least as meta-theorems for small sets, or for sets which can be constructed as closed terms in ZFC. [I.e., the term model of ZFC is a part of the term model of NACT+NFUM().] Probably, the axioms of ZFC restricted to small sets are also valid in NACT+NFUM(). Therefore NACT+NFUM() can also be considered as an enlargement of ZF''. (ZF''[21] is a slight modification and restriction of ZF [with the axioms of separation and replacement restricted to small (here: well-founded) sets and an ω-axiom instead of the axiom of infinity.] In ZF(4) also an axiom for the Cartesian product exists and the (AC) is also restricted.)[22]

---

[20] See DePauli [2007b] "Paradigmen-Wechsel in der Mengenlehre" in EUROPOLIS4.

(10)  Free Parameters in the Class Operators

Let us now consider the system NACT☼W (= the naked sun Wittgenstein). After the explanation of the construction of list2 and generation-algorithm GL2 and of axiom-list2 we can now proceed to list3, GL3 and axiom-list3. First we have to modify algorithm GL2 in such a way that also free parameters b1, b2, b3, etc. can be included in the formulas. But since the generation of formulas in list3 is synchronized with the length of the formulas (as in list1 or 2 as well), the parameters can only be introduced in a way similar to the bounded variables (i.e., in every circle only one new one).

Next, we have to solve the question of what it means to prove the main task for a formula $A_i(x, b_1, b_2, \ldots b_n)$ ? I.e. for NACT☼W:

? Set($\{|x: A_i(x, b_1, b_2, \ldots b_n)|\}$) $\vdash_{\overline{CT\&\ ES}}$ falsum ?

Let us first consider this question for a single parameter b. It is equivalent to show:

? Set($\{|x: A_i(x, b)|\}$) $\vdash_{\overline{CT\&\ ES}}$ falsum ?

[E.g., find a set b such that Set($\{|x: x \notin b|\}$) $\vdash_{\overline{CT\&\ ES}}$ falsum.] This means in practice that we have to prove A(x, b) [with the limitation length(A) ≤ a given n] for all already proven sets s with length(s) ≤ length(A) [where length(s) means the length of the set-constituting formulas]. Therefore the formula A(x, s) can only become twice as long as A(x, b) and so this proof is limited too.

But this does not give us any certainty, because we should prove this for all sets s. [And furthermore, in our computational set theory all sets are constructed as terms. So we assume implicitly the term model that every set can be constructed by a term[23]. But this may be a good thing, because it simplifies set theory.]

Therefore, as long as we cannot prove the invalidity of the main task for formulas with parameters for all substitution instances, we can in practice only prove it for instances with a limited length and assume only a provisional sethood of such a class. [Because if further substitution instances affirm the main goal, it will turn out that such classes with parameters

can become properclasses in a revision. E.g., if we investigate the "complement of the Russell class".]

Since this is only an Gedankenexperiment, let us agree to construct axiom-list3 in the following way: the parameters in formula $A_i(x, b_1, b_2, \ldots b_n)$ have to be assigned with all combinations of already proven sets. This has to be done again and again after the proof of existence of new sets. The formula $A_i$ should not grow over a given limit [i.e. the longest derived formula without parameters of list3, which is here identical with list2]. Thus, the assignment with all combinations of a fixed length cannot force a combinatorial explosion. This way we can arrange that the assigned formula has been already produced by list3 or list2, respectively. But the price for this trick is that the axioms with parameters in fact remain forever undecided (if we cannot find some meta-mathematical proof to the contrary).

Since NACT☼W allows parameters, it is possible to prove Cantor's theorem. Furthermore, the complement to every concretely given set exists, except to the complements of the properclasses.

(11)   The System NACT☼

Let us now simplify NACT☼W and cross out ES under the derivation sign of the main task. This changes the main task to

? Set({|x: $A_i(x, b_1, b_2, \ldots b_n)$|}) $\vdash_{CT}$ falsum ?

So the main task of NACT☼ is slightly easier than that of NACT☼W, but the problem with the substitution of all combinations of existing sets still remains. Therefore we do not have in fact any benefit compared with NACT☼.

NACT☼ is more or less NACT* with parameters in the set-constituting class operators. It would also be possible to take NACT ❸ with parameters, but we cannot find a formal system as an approximation to NACT ❸ . Therefore we skip this possibility and construct NACT☼ [on the basis of list2 and GL2] as extension of NACT*. The study of the extension of NACT ❸ would not be very profitable.

Furthermore: if we want to combine NACT☼ with NFUM in order to construct an approximation similar to NACT+NFUM() [to NACT*], we must be careful, since NFUM has a complement. Therefore, we may not use (J3) in the same way as before. The first solution we suggest is to use ZFC4 (of chapter 9) instead of NFUM or (J2) & (J3) & ZFC4& (GCH). The second is to use a 2-level system where (J3) is only allowed to be applied to the finished system NFUM [i.e., no axiom of NFUM can be used after (J3)].

The 1st solution ZFC4 & (GCH) & (J2) & (J3) does not support our initial goal of constructing a system with the "most possible number of sets". The 2nd solution is not a standard formal system and leaves a lot of questions unanswered. But it is necessary to find a system with parameters in those class operators should establish a set. [In general classes which need not be sets, parameters are of course allowed.]

(12) The System NACT+4, NACT & ZFC4+(GCH)

NACT+4 has the following axioms:[24]

1:    Class-Operator
2:    Church Schema
3:    Extensionality
4:    a weak (AC)
5a:   Disjunctive sets axiom: Set(X) or Set($\neg$X).
6c:   Ordering axiom: Small(X) $\rightarrow$ Set(X) and Set ($\neg$X).
      Here Small can be a suitable properly like Slim, Well-Founded; Cantorian, etc.
      We use Slim(X) : = |X| < |$\neg$X|.

To these 6 basic axioms we add the 4 "ZF-axioms":

(7#)   Small ($\omega_0$),
(8#)   Small (X) $\rightarrow$ Small(Power(X)),
(9#)   Small (X) [& $\forall y \in$ X: Small(y)] $\rightarrow$ Small ($\bigcup$X), where "$\bigcup$" is the large union.
(10#)  Small(X) & Function (F) $\rightarrow$ Small (F[X]), where F[X] is the image of X under F.

Since we always assume the validity of the Class-Theoretical frame, axioms 5a and 6c are (J3) and (J2), and ZF4 are the 4 ZF-axioms (7#) to (10#), we can call this system also ZFC4 & (J2) & (J3). If we add (GCH) we can also call this system ZFC4 & (GCH) & (J2) & (J3). (J2) generates here additional sets which do not exist in ZFC4 -- e.g., $\mathbb{V}$ = $\neg$0. But Tarski's unaccessible number $\tau_0$ can also not be constructed as set in ZFC4, while it is probably a set in NACT & ZFC4+(GCH). Since its complement contains $\neg$Ru [= the complement of the Russell class] by (J3) and all its sets [containing themselves], $\tau_0$ is of smaller power than $\neg\tau_0$ and therefore a set concerning (J2).

But what does it mean to say that "X is of smaller power than Y"? Answer: There is no 1-1 function from X onto Y. These functions need not be sets here. [It is enough if they are classes or formulas with the functional property.] This is no handicap, but it weakens the system a little bit. It is disadvantageous if we should want to show that some specific class is medial, or that the medial classes are equally mighty to the universal set. But we do not even know if any medial class is a set in NACT & ZFC4+(GCH).

All in all, this system is unsatisfying to the imagination, and we should therefore proceed to the next system.

(13)   The System NACT & NFUM

[Let us call this system also NFUM & (J2) / (J3).] It is not clear if (J2) is derivable from NFUM [and therefore redundant]. Therefore we want to mention it in any case. The slash between (J2) and (J3) means that NFUM & (J2) has first to be worked out. After that (J3) can be used. But after using (J3), no axiom of NFUM can be applied again. Especially not the complement, because by (J3), $\neg$Ru is a set. And if the complement of $\neg$Ru were to be a set,

---

we would have the Russell Paradox. So NFUM[&(J2)]/(J3) is a 2-level system or also a 2-sorted formal system where (J3) can be only applied to the $1^{st}$ sort of entities or variables of NFUM&(J2). The result of (J3) are the $2^{nd}$ sort entities.

In this system medial sets also exist, e.g. {x: 0 ∉ x}, which is stratified and of equal power to its complement. The function between them is also a set. The correspondence between {x: 0 ∉ x} and 𝕎 is probably also a set. This yields the somehow counter-intuitive result that many medial classes are equipollent to the universe. But that is no contradiction.

After having understood this system, we can ask ourselves under what condition we can apply again an axiom of NFUM? E.g.: it will certainly be allowed to form {⌐Ru} [= the singleton of the complement of the Russell class]. Therefore we can add axioms (X1), (X2), . . . (Xn) to the system, yielding NFUM&(J2)/(J3)&<(Xi)>. These series of axioms <(Xi)> are of course stratified, but their set-constituting formula $A_i$ should not contain a subformula B from which the existence of the complement can be derived. This is probably (but not certainly) an undecidable property, and we have again the same problems we already had in the solution of the main task of NACT ☻, *, and ☼.

Nevertheless, its already enough to have NFUM&(J2)/(J3). The additional axioms <(Xi)> are not so important. [But they can generate sets which can possibly be joint consistently to this system.] Therefore we consider the originally proposed goal of constructing a "greatest possible maximal set theory" as solved (for the next 10 years).

(14) The Systems NACT&ZFC4+(GCH) and NACT[&ZFC4]+(FCA)
If we renounce our original goal to construct the "greatest possible maximal set theory" and are satisfied with a "suitably large set theory" [which is therefore larger than ZFC] we need not to work with NFUM and consider now the formal system (J2) & (J3) &ZFC4 & (GCH) [with parameters].

The 2 principles (J2) and (J3) can blow up a universe only as large as ω. But it is not certain that we can prove Cantor's Theorem with them. Also ZFC4-closed does not help. Neither does the use of General Continuum Hypothesis[25]: "$|X| = \aleph_\alpha$ ➔ $|P(X)| = \aleph_{\alpha+1}$ create a larger universe. [It may turn out that the Axiom of Choice is only valid for small sets and we have to use (AC-small) instead of (AC).]

We are here in the same situation as with NACT ☻ and NACT*: Instead of (GCH) we can also use the "Finite or Countable Anti-Thesis" (FCA): $|X| \leq \omega_0$.[26] But, in this case, it may possibly be better to also restrict (J2) and (J3) to parameter-free formulas and/or (FCA) and

---

[25] See Gödel [1956]: What is Cantor's Continuum Problem? And Gödel [1938].

[26] Or another anti-thesis such as e.g: ($|X| = \omega_0 => |P(X)| = \omega_1$) & ∀X: $|X|$ = finite or $\omega_0$ or $\omega_1$). See also Randall Holmes Pocket Theory.

(AC) to small sets. We can use ZFC4-closed or not, as we like[27]. In this connection many problems are still left unsolved, and much research has yet to be done.

   (15) SUMMARY

We have tried in this article to find a new paradigm for set theory. All classes which are possible to considered consistently as sets (in coexistence with the common sets of ZFC) should be sets and can also be established as sets. NACT☼ and its approximation NACT&NFUM are also models for the practical work of mathematicians [who are working in practice with Naïve Set Theory based on Non-Monotonic Reasoning]. It shows that there is no danger of encountering inconsistency. Furthermore: the computational theorem proving[28] with these systems is much easier than with ZFC.

   In addition, there are a lot of non-trivial technical problems open that seem interesting. But of course, the reader has first to study NFUM. The book of Randall Holmes is out of print, but you can load it down from his homepage. I guarantee that it is really worth reading! Presently NFUM is the basis of all "maximal set theories" which should work like a formal system. [This is the result of the author's many-years' work in this field.]

   But for the beginning, we can also try to avoid the use of NFUM and use the philosophical principles (J2) "Small classes and their complements are sets" and (J3) "Every class or its complement is a set, or both" together with ZFC4 and the General Continuum Hypothesis. This also yields a maximal set theory, but not the greatest possible one.*

   (16) Literature

---

[27] An alternative would be to use ZFC4 and restrict the use of the axiomy of separation and replacement to definite formulas obeying the rules of the theory of definition. See DePauli-Schimanovich [2008c].

[28] See "Georg Gottlob, Alexander Leitsch, Werner Schimanovich [1986]: Automatisches Beweisen (Methoden und Implementation)" in the literature DePauli [2005a] "The Vienna Theorem Prover" in EUROPOLIS5.

* Acknowledgements:

I want to thank Matthias Baaz, Martin Goldstern and Georg Gottlob for several advices for this article.